\documentclass[12pt]{amsart}

\usepackage{amsmath, amssymb}

\newtheorem{theorem}{Theorem}

\newtheorem{lemma}[theorem]{Lemma}

\newtheorem{definition}{Definition}

\def\qedf{\hfill $\Box$}

\title{{\bf Weighted Regularity Lemma with Applications}}

\author[Csaba]{B\'ela Csaba}
\address{Department of Mathematics, Western Kentucky University}
\email{bela.csaba@wku.edu}

\author[Pluh\'ar]{Andr\'as Pluh\'ar}
\address{University of Szeged, Department of Computer Science, \'Arp\'ad t\'er 2., Szeged, H-6720}
\email{pluhar@inf.u-szeged.hu}

\thanks{The authors were partially supported by the grants OTKA T049398 and OTKA K76099. 
The first author was partially supported by a New Faculty Scholarship Grant of the Office of
Sponsored Programs at Western Kentucky University.}
\keywords{Regularity lemma, sparse graphs, quasi-random.}

\subjclass {}

\date{}

\begin{document}

\begin{abstract} 

We prove an extension of the Regularity Lemma with vertex and edge weights which can be applied for a large class of graphs.
The applications involve random graphs and a weighted version of the Erd\H{o}s-Stone theorem. We also provide means to handle the otherwise uncontrolled exceptional set.
\end{abstract}

\maketitle

\section{Introduction}

Let $G=G(A,B)$ be a bipartite graph. For $X, Y \subset A \cup B$ let $e_G(X, Y)$ denote the number of edges with one endpoint in
$X$ and the other in $Y.$
We say that the $(A, B)$-pair is $\varepsilon$-regular if 

$$\left| \frac{e_G(A', B')}{|A'| |B'|}-\frac{e_G(A, B)}{|A|  |B|} \right|< \varepsilon$$
for every $A \subset A,$ $|A'|>\varepsilon |A|$ and $B' \subset B,$ $|B'| > \varepsilon |B|.$

This definition plays a crucial role in the celebrated Regularity Lemma of Szemer\'edi, see~\cite{SZ75, SZ76}. 
Szemer\'edi invented
the lemma for proving his famous theorem on arithmetic progression of dense subsets of the natural numbers.
The Regularity Lemma is a very powerful tool when applied to a dense graph. It has found lots of applications
in several areas of mathematics and computer science, for applications in graph theory see e.g.~\cite{KS93}.  
However, it does not tell us anything useful when 
applied for a sparse graph (i.e., a graph on $n$ vertices having $o(n^2)$ edges).

There has been significant interest to find widely applicable versions for sparse graphs. This turns out to be
a very hard task. Kohayakawa~\cite{K} proved a sparse regularity lemma, and with R\"odl and \L uczak~\cite{KLR} they applied
it for finding arithmetic progressions of length 3 in dense subsets of a random set. In their sparse regularity
lemma dense graphs are substituted by dense subgraphs of a random (or quasi-random) graph. Naturally, a new definition
of $\varepsilon$-regularity was needed, below we formulate a slightly different version from theirs:

Let $F(A,B)$ and $G(A,B)$ be two bipartite graphs such that $F \subset G.$ We say
that the $(A,B)$-pair is $\varepsilon$-regular in $F$ relative to $G$ if
$$\left| \frac{e_F(A', B')}{e_G(A', B')}-\frac{e_F(A, B)}{e_G(A, B)} \right|< \varepsilon$$
for every $A' \subset A, B' \subset B$ and $|A'| > \varepsilon |A|$, $|B'| > \varepsilon |B|$.
It is easy to see that the above is a generalization of $\varepsilon$-regularity, in the original definition 
the role of $G$ is played by the complete bipartite graph $K_{A, B}.$ In this more general definition $F$ can 
be a rather sparse graph, it only has to be dense {\it relative to $G$}, that is, $e(F)/e(G)$ should be a constant.

In this paper we further generalize the notion of quasi-randomness and $\varepsilon$-regularity by introducing 
{\em weighted regularity} using vertex and edge weights. This enables us to prove a more general, and perhaps 
more applicable regularity lemma. 
Let us remark, that a new notion, {\em volume regularity} was recently used by Alon et al.~\cite{A}. As we will
see later, volume regularity is a special case of weighted regularity. 
Still, these approaches are not fully comparable, since they use different premises for finding the regular partition of a graph. In this
paper we discuss the Regularity Lemma for dense subgraphs of {\em weighted quasi-random} graphs, while Alon et al. find
a regular partition for so called $(C, \eta)$-{\em bounded} graphs. 

The basic tool is the Strong Structure Theorem of Tao~\cite{TT07}, where he simplifies the proof of the original
Regularity lemma itself and gives new insights, too. Following his lines became technically feasible to extend
regularity to the case when both the edges and the vertices of a graph are weighted. (Note, that the measures 
are in close connection with each other.) 
We remark that similar ideas might be used to find a regularity lemma for certain sparse hypergraphs as well.

The structure of the paper is as follows. First we discuss weighted quasi-randomness and weighted $\varepsilon$-regularity in the second
section. In the third section we prove the new version of the regularity lemma. Finally, we show some applications in the fourth section, in particular,
we prove a weighted version of the Erd\H{o}s-Stone theorem.

\section{Basic definitions and tools}

Throughout the paper we apply the relation ``$\ll$'':  $a\ll b$ is $a$ is sufficiently smaller than $b.$ This notation is widely applied in papers 
using the Regularity Lemma, and simplifies our notation, too.

Let $\beta >0$ and $G=(V,E)$ be a graph on $n$ vertices. Set $\delta_G= {e(G)}/{\binom{n}{2}},$ this is the 
{\it density} of $G.$ We define the density of the $A, B$ pair of subsets of $V(G)$ by 
$\delta_G(A, B)={e_G(A, B)}/(|A| |B|).$ We say that $G$ is 
$\beta$-quasi-random if it has the following property: 
If $A, B \subset V(G)$ such that $A\cap B=\emptyset$ and $|A|, |B| > \beta n$ then
$$\left|\delta_G-\delta_G(A, B)\right|<\beta \delta_G.$$
That is, the edges of $G$ are distributed ``randomly.'' In order to formulate our 
regularity lemma we have to define quasi-randomness in a more general way, 
that admits weights on vertices and edges.

For a function $w: S \rightarrow \mathbb R^+$  
and $A \subset S$, $w(A)$ is defined by the usual way, that is, 
$w(A)=\sum_{x \in A}w(x).$ 
We shall also use the indicator function of the edge set of a graph $H.$ 
${\mathbf 1}_H$: $\binom{V(H)}{2} \rightarrow \{0, 1\}$ and 
${\mathbf 1}_H(x,y)=1$ iff $xy \in E(H).$

We define the weighted quasi-randomness of a graph $G=(V, E)$ with given weight-functions 
$\mu: V \rightarrow \mathbb R^+$
and $\rho: \binom{V}{2}\rightarrow \mathbb R^+.$ For $A, B \subset V$ let 
$$\rho_G(A, B):=\sum_{u \in A, v \in B}{\mathbf 1}_G(u, v)\rho(u, v).$$  
In particular, $\rho_G(u,v)={\mathbf 1}_G(u,v)\rho(u,v)$ for $u, v \in V.$
Observe, that the function $\mu$ is an analogon of the vertex counting function on a set, while
the function $\rho$ counts the edges in the unweighted case.

\begin{definition} \label{quasirandom}
A graph $G=(V, E)$ is weighted $\beta$-quasi-random with weight-functions $\mu$ and $\rho$ if
for any $A, B \subset V(G)$ such that $A \cap B =\emptyset$ and 
$\mu(A) \geq \beta \mu (V)$, $\mu(B) \geq \beta \mu (V)$ we have

$$\left|\frac{\rho_G(A, B)}{\mu(A)\mu(B)}- 
\frac{\rho_G(V, V)}{\mu(V)\mu(V)}\right| < \beta.$$
\end{definition}

Observe that choosing $\mu \equiv 1$ and $\rho \equiv 1/\delta_G$ gives back the 
first definition of quasi-randomness. There is another, weaker notion of quasi-randomness,
which we will also use. 

\begin{definition} \label{D-quasirandom}
Let $D>1$ be an absolute constant. A graph $G=(V, E)$ is weighted $(D, \beta)$-quasi-random with weight-functions $\mu$ and $\rho$ if
for any $A, B \subset V(G)$ such that $A \cap B =\emptyset$ and 
$\mu(A) \geq \beta \mu (V)$, $\mu(B) \geq \beta \mu (V)$ we have

$$\frac{1}{D}\frac{\rho_G(V, V)}{\mu(V)\mu(V)} \le \frac{\rho_G(A, B)}{\mu(A)\mu(B)} \le 
D \frac{\rho_G(V, V)}{\mu(V)\mu(V)}.$$
\end{definition}

Clearly, if a graph is $\beta$-quasi-random and $D > \max\{1+\beta/y, 1+\beta/(y-\beta)\}$,
then it is $(D, \beta)$-quasi-random, where $y:=\rho_G(V, V)/\mu(V)^2$.
Now we need to describe the weighted version of relative regularity. 

\begin{definition}
Let $G$ and $F$ be graphs such that $F \subset G$ and $\mu, \rho$ weight 
functions defined as above. 
For $A, B \subset V(G)$ and $A \cap B =\emptyset$ the 
pair $(A, B)$ in $F$ is $(\mu, \rho)$-weighted $\epsilon$-regular relative to $G$, 
or briefly weighted $\epsilon$-regular, if
$$\left|\frac{\rho_F(A', B')}{\mu(A')\mu(B')}- 
\frac{\rho_F(A, B)}{\mu(A)\mu(B)}\right| < \epsilon$$
for every $A' \subset A$ and $B' \subset B$ provided that 
$\mu(A') \geq \epsilon \mu(A), \mu(B') \geq \epsilon \mu(B)$.
Here 
$$\rho_F(A, B)=\sum_{u \in A, v \in B}{\mathbf 1}_F(u, v)\rho(u, v).$$ 
\end{definition}

\noindent {\bf Remarks.} Note that weighted $\epsilon$-regularity is nothing but 
the well-known $\epsilon$-regularity when $G=K_{A, B}$ and $\mu \equiv 1$ and $\rho$
is chosen to be identically the reciprocal of the density of $G$ as before.
Since ${\mathbf 1}_F(u, v) \leq {\mathbf 1}_G(u, v) \leq {\mathbf 1}(u, v)$ we also have 
$\rho_F(A, B) \leq \rho_G(A, B) \leq \rho(A, B)$. Hence, the first inequality of the definition
does not refer to $G$ explicitly, but contains information on it.

Next we define weighted regular partitions.
\begin{definition} \label{wregular} 
Let $G=(V, E)$ and $F\subset G$ be graphs, and $\mu$ and $\rho$ weight functions. $F$ has a 
weighted $\epsilon$-regular partition relative to $G$ if its vertex set $V$ can be partitioned 
into $\ell +1$ clusters $W_0, W_1, \dots, W_{\ell}$ such that

\begin{itemize} 
\item $\mu(W_0) \leq \epsilon \mu(V),$
\item $|\mu(W_i)-\mu(W_{j})| \le \max_{x \in V}\{\mu(x)\}$ for every $1 \le i, j \le \ell,$
\item all but at most $\epsilon \ell^2$ of the pairs $(W_i, W_j)$ for
$1 \leq i < j \leq \ell$ are weighted $\epsilon$-regular in $F$ relative to $G.$
\end{itemize}

\end{definition}

In order to show our main result we will use the Strong Structure Theorem of Tao, 
that allows a short exposition. In fact we
will closely follow his proof for the Regularity Lemma as discussed in~\cite{TT07}. 
First we have to introduce some definitions.

\smallskip

Let $H$ be a real finite-dimensional Hilbert space, and let $S$ be a set of basic functions
or basic structured vectors of $H$ of norm at most 1. The function $g \in H$ is 
$(M,K)$-structured with the positive integers $M, K$ if one has a decomposition
$$g=\sum_{1 \le i \le M}{c_i s_i}$$
with $s_i \in S$ and $c_i \in [-K, K]$ for $1 \le i \le M.$ We say that $g$ is 
$\beta$-pseudorandom for some $\beta >0$ if $|\langle g, s \rangle| \le \beta$ for 
all $s \in S$. Then we have the following

\begin{theorem}[Strong Structure Theorem - T. Tao] \label{structure}
Let $H$ and $S$ be as above, let $\varepsilon >0$, and let $J: {\rm {\bf Z}}^+ \rightarrow {\mathbb R}^+$ 
be an arbitrary function. Let $f \in H$ be such that $\|f\|_H\le 1.$ Then we can find an integer 
$M=M_{J,\varepsilon}$ and a decomposition $f=f_{str}+f_{psd}+f_{err}$ where
(i) $f_{str}$ is $(M,M)$-structured, (ii) $f_{psd}$ is $1/J(M)$-pseudorandom, and (iii) 
$\|f_{err}\|_H\le \varepsilon.$ 
\end{theorem}

\section{Weighted Regularity Lemma relative to a quasi-random graph $G$}

First we define the Hilbert space $H$, and $S$. We generalize Example 2.3 of~\cite{TT07} to weighted
graphs.
Let $G=(V, E)$ be a $\beta$-quasi-random graph on $n$ vertices with weight functions $\mu$ and $\rho.$ 
Let $H$ be the $\binom{n}{2}$-dimensional space of functions
$g: \binom{V}{2} \rightarrow {\mathbb R},$ endowed with the inner product
$$\langle g, h \rangle = \frac{1}{\binom{n}{2}}\sum_{(u,v) \in \binom{V}{2}}g(u,v)h(u,v)\rho_G(u,v).$$

It is useful to normalize the vertex and edge weight functions, we assume that 
$\mu(V)=n$ and $\langle 1, 1\rangle =1.$ We also assume, that $\mu(v)=o(|V|)$ for every $v \in V.$
Observe, that if $F \subset G$ then $\|{\mathbf 1}_F\|\le 1.$ 
We let $S$ to be the collection of 0,1-valued functions $\gamma_{A, B}$ for 
$A, B \subset V(G)$, $A \cap B=\emptyset$, where $\gamma_{A,B}(u, v)=1$ if 
and only if $u \in A$ and $v \in B.$ We have the following

\begin{theorem}[Weighted Regularity Lemma] \label{grl}
Let $D>1$ and $\beta, \varepsilon \in (0,1),$ such that $0 < \beta \ll \varepsilon \ll 1/D$ and let $L \ge 1.$ 
If $G=(V, E)$ is a weighted $(D,\beta)$-quasi-random graph on $n$ vertices with $n$
sufficiently large depending on $\varepsilon$ and $L$ and $F \subset G$, 
then $F$ admits a weighted  $\varepsilon$-regular partition relative to $G$
into the partition sets $W_0, W_1,  \ldots ,W_{\ell}$ such that 
$L \le \ell \le C_{\varepsilon, L}$ for some constant $C_{\varepsilon, L}.$
\end{theorem}

\noindent {\bf Proof:}
Let us apply Theorem~\ref{structure} to the function ${\mathbf 1}_F$ with parameters 
$\eta$ and function $J$ to be chosen later.
We get the decomposition $${\mathbf 1}_F=f_{str}+f_{psd}+f_{err},$$
where $f_{str}$ is $(M,M)$-structured, $f_{psd}$ is $1/J(M)$-pseudorandom, and $\|f_{err}\| \le \eta$ with
$M=M_{J, \eta}=M_{J, \varepsilon}.$

The function $f_{str}$ is the combination of at most $M$ basic functions:
$$f_{str}=\sum_{1 \le k \le M}\alpha_k \gamma_{{\mathcal A}_k, {\mathcal B}_k}$$
where ${\mathcal A}_k, {\mathcal B}_k$ are subsets of $V$ and $\gamma_{{\mathcal A}_k, {\mathcal B}_k}$ 
agrees with the indicator function of the edges of $G$ in between ${\mathcal A}_k$ and ${\mathcal B}_k.$ 
Any $({\mathcal A}_k,{\mathcal B}_k)$ pair partitions $V$ into at most 4 subsets. Overall we get 
a partitioning of $V$ into at most $4^M$ subsets, we will refer to them as {\em atoms.} Divide every 
atom into subsets of total vertex weight $\frac{\varepsilon n}{ L+4^M}$, except possibly one smaller 
subset. The small subsets will be put into $W_0$, the others give $W_1, W_2, \ldots, W_{\ell}$, with 
$\ell=\frac{L+4^M}{\varepsilon}$. We refer to the sets $W_i$ for $i=1, \dots \ell$ as {\it clusters}. 
If $n$ is sufficiently 
large then this partitioning is non-trivial. From the construction it follows that each $W_i$ is entirely 
contained within an atom. It is also clear that $\mu(W_0) \le \varepsilon n$ and 
$\mu(W_i)\approx m=\Theta(\frac{n}{\ell})$ for every $1 \le i \le \ell.$

We have that 
$$\|f_{err}\|^2=\frac{1}{\binom{n}{2}}\sum_{(u, v) \in \binom{V}{2}}|f_{err}(u,v)|^2\rho_G(u,v) \le \eta^2.$$
From this and the normalization of $\rho$ it follows that  
$$\frac{1}{\binom{\ell}{2}} \sum_{1 \le i<j \le \ell}\frac{1}{\rho_G(W_i, W_j)}\sum_{u\in W_i, v \in W_j}|f_{err}(u,v)|^2\rho_G(u,v)=O(\eta^2).$$
Clearly,
$$\frac{1}{ \rho_G(W_i, W_j)}\sum_{u \in W_i, v \in W_j}|f_{err}(u,v)|^2\rho_G(u,v)=O(\eta) $$
for all but at most $O(\eta \ell^2)$ pairs $(i,j).$ If the above is satisfied for a pair $(i,j)$ then we call it a {\it good pair.} 
We will apply the Cauchy-Schwarz inequality. For that let 
$a(u,v)=|f_{err}(u,v)|\sqrt{\rho_G(u,v)}$ and $b(u,v)=\sqrt{\rho_G(u,v)},$ then 

$$\frac{\sum_{u \in W_i,v\in W_j}a(u,v)b(u,v)}{\sqrt{\sum_{u\in W_i, v\in W_j}b^2(u,v)}}\le \sqrt{\sum_{u\in W_i, v\in W_j}a^2(u,v)}.$$
Since
$$\sqrt{\sum_{u\in W_i, v\in W_j}a^2(u,v)}=O(\sqrt{\eta})\sqrt{\rho_G(W_i, W_j)},$$
we get that $$\frac{1}{\rho_G(W_i, W_j)}\sum_{u\in W_i, v\in W_j}|f_{err}(u,v)|\rho_G(u,v) =O(\sqrt{\eta})$$
if $(i,j)$ is a good pair.

Assume that $(i,j)$ is a good pair. From the pseudorandomness of $f_{psd}$ we have that 
$$|\langle f_{psd}, \gamma_{A, B} \rangle| = \frac{1}{\binom{n}{2}}\left|\sum_{u\in A, v\in B}f_{psd}(u,v)\rho_G(u,v)\right| \le \frac{1}{ J(M)}$$
for every $A \subset W_i$ and $B \subset W_j.$ 

\medskip

We will show that every good pair is weighted $\varepsilon$-regular in $F$ relative to $G.$ 
Let $(i, j)$ be a good pair, and assume that $A \subset W_i,$ $\mu(A)>\varepsilon \mu(W_i)$ and $B \subset W_j,$ $\mu(B)>\varepsilon \mu(W_j).$
To show that $(W_i,W_j)$ is weighted $\varepsilon$-regular, it is sufficient to show that

$$\left| \frac{\rho_F(A, B)}{\mu(A)\mu(B)} - \frac{\rho_F(W_i, W_j)}{\mu(W_i)\mu(W_j)} \right|<\varepsilon.$$
Recall that $$\rho_F(A,B)=\sum_{u \in A, v \in B}{\mathbf 1}_F(u,v)\rho(u,v)=\sum_{u \in A, v \in B}{\mathbf 1}_F(u,v)\rho_G(u,v),$$
since $F \subset G.$

Substituting $f_{str}+f_{psd}+f_{err}$ for ${\mathbf 1}_F$ it is sufficient to verify the following inequalities.

$$(1) \left| \frac{\sum_{u \in A, v\in B}f_{str}(u,v)\rho_G(u,v)}{\mu(A)\mu(B)}-\frac{\sum_{u \in W_i, v\in W_j}f_{str}(u,v)\rho_G(u,v)}{\mu(W_i)\mu(W_j)}\right|<\varepsilon/3,$$
$$(2) \left| \frac{\sum_{u \in A, v\in B}f_{psd}(u,v)\rho_G(u,v)}{\mu(A)\mu(B)}-\frac{\sum_{u \in W_i, v\in W_j}f_{psd}(u,v)\rho_G(u,v)}{\mu(W_i)\mu(W_j)}\right|<\varepsilon/3$$
and 
$$(3) \left| \frac{\sum_{u \in A, v\in B}f_{err}(u,v)\rho_G(u,v)}{\mu(A)\mu(B)}-\frac{\sum_{u \in W_i, v\in W_j}f_{err}(u,v)\rho_G(u,v)}{\mu(W_i)\mu(W_j)}\right|<\varepsilon/3.$$

\medskip

For proving (1) recall that $f_{str}$ is constant on $W_i \times W_j$ and $(M,M)$-structured. 
Since the $\gamma_{X, Y}$ basic functions are $0, 1$-valued, we get, that $|f_{str}|\le M^2.$ 
Moreover, $G$ is $(D, \beta)$-quasi-random, where $0<\beta \ll \varepsilon.$ Therefore, 
$(1) \le D M^2\beta < \varepsilon/3$, since $\beta \ll \epsilon$.

The proof of (2) goes as follows. The first term is
$$\left| \frac{\sum_{u \in A, v\in B}f_{psd}(u,v)\rho_G(u,v)}{\mu(A)\mu(B)} \right|=
\binom{n}{2}|\langle f_{psd}, \gamma_{A, B} \rangle | \le \frac{\binom{n}{2}}{J(M)\mu(A)\mu(B)}$$

the second is 

$$\left| \frac{\sum_{u \in W_i, v\in W_j}f_{psd}(u,v)\rho_G(u,v)}{\mu(W_i)\mu(W_j)} \right| =
\binom{n}{2}|\langle f_{psd}, \gamma_{W_i, W_j} \rangle |
\le \frac{\binom{n}{2}}{J(M)\mu(W_i)\mu(W_j)}.$$
Noting that $\mu(W_k) = \Theta(n/\ell)$ for $k \ge 1$ we get that the sum of the above terms is at most 
$$\frac{\ell^2}{2J(M)}\left(1+\frac{1}{\varepsilon^2}\right)< \frac{\varepsilon}{3},$$
if $J(M) \gg \frac{\ell^2}{\varepsilon^3}.$

For $(3)$ first notice that it is upper bounded by 
$$O(\sqrt{\eta})\left( \frac{\rho_G(W_i, W_j)}{\mu(W_i)\mu(W_j)}+\frac{\rho_G(W_i, W_j)}{\mu(A)\mu(B)} \right)\le 
O(\sqrt{\eta})\frac{\rho_G(W_i, W_j)}{\varepsilon^2\mu(W_i)\mu(W_j)}.$$ 
We also have that
$$\frac{\rho_G(W_i, W_j)}{\mu(W_i)\mu(W_j)}=O(1)$$
by the normalization of $\mu$ and $\rho$ and from the fact that $G$ is quasi-random. From this it is easy to see 
that if $\eta \ll \varepsilon^6$ then $(3)$ is at most $\varepsilon/3.$
This finishes the proof of the theorem. \hfill \qedf

\section{Applications}

In this section we consider a few applications of our main result. First we prove that a random graph with 
widely differing edge probabilities 
is quasi-random, if none of the edge probabilities are too small. In this case the vertex weights will all be one, 
but edges will have different weights. Then we show examples where vertices have different weights. We will
consider the relation of weighted regularity and volume regularity. We define the `natural weighting' of $K_n,$
and prove a weighted version of the Erd\H{o}s-Stone theorem for this weighting. Finally, we show how to
partially control the non-exceptional set by natural weightings. 

\subsection{Quasi-randomness in the $G(n, p_{i j})$ model}

In this section we will prove that random graphs of the $G(n, p_{ij})$ model are quasi-random in the strong sense with high probability. 
A special case of this model is the well-known $G(n,p)$ model for random graphs. A Regularity Lemma for this case was 
first applied by Kohayakawa, \L uczak and  R\"odl in~\cite{KLR}. They studied $G(n, p)$ for 
$p=c/\sqrt{n}$ in order to find arithmetic progressions of length three 
in dense subsets of random subsets of $[N]$.

The $G(n, p_{i j})$ model was first considered by Bollob\'as \cite{Bollobas}.
Recently it was also studied by Chung and Lu \cite{CL}.

In this model one takes $n$ vertices, and draws an edge
between the vertices $x_i$ and $x_j$ with probability $p_{ij}$, randomly 
and independently of each other. Note that if $p_{ij} \equiv p$, then 
we get back the well-known  $G(n, p)$ model. It is a straightforward application of the Chernoff bound
that a random graph $G \in G(n, p)$ is quasi-random with high probability. However, the case of $G(n, p_{i, j})$
is somewhat harder.

\begin{lemma} \label{pij} 
Let $\beta>0.$ There exists a $K=K(\beta)$ such that if $G \in G(n, p_{i j})$ and 
$p_{i j} \ge K/n$ for every $i$ and $j$, then $G$ is weighted $\beta$-quasi-random 
with probability at least $1-2^{-n}$ if $n$ is sufficiently large. 
\end{lemma}

\noindent {\bf Proof.} First of all let $\mu \equiv 1,$ and let $\rho (i,j)=1/p_{i,j}.$ 
Set $K=4800/\beta^6.$
Let $p_0=K/n$, and let $p_k=e^{k}p_0$ for $1 \le k \le \log n.$ 
Let $A$ and $B$ be a pair of disjoint sets, both of size at least $\beta n.$
We partition the pairs $(u,v)$, where $u \in A$ and $v \in B,$ into $O(\log n)$ 
disjoint sets $H_1, H_2, \ldots, H_l:$ 
if $p_{k}\le p_{uv} <p_{k+1}$ then $(u,v)$ will belong to $H_k.$ 
Let $a_k= \frac{\beta^3}{10} \sqrt{e}^{k}Kn.$ We will denote $|H_k|$ by $m_k.$

We will prove that the following inequality holds with probability at least $1-2^{-3n}:$ 

$$\left|\sum_{u \in A, v\in B} \frac{X_{uv}}{p_{uv}|A||B|}-1\right|<\beta/2,$$
where $X_{uv}$ is a random variable which is 1 if $uv \in E(G),$ otherwise it is 0. This
implies the quasi-randomness of $G$ since there are less than $2^{2n}$ pairs of disjoint subsets of $V(G).$ Observe that 
$$ \sum_{u \in A, v\in B} \frac{{\mathbb E}X_{uv}}{p_{uv} |A||B|}=1.$$

Applying the large deviation inequalities $A.1.11$ and $A.1.13$ from~\cite{A-S}, 
we are able to bound the number of edges in between $A$ and $B$ for the edges of $H_k$ in case $m_k$
is sufficiently large as follows. According to
$A.1.11$ we have that 
$${\rm Pr}\left(\sum_{(u,v)\in H_k}(X_{uv}-{\mathbb E}X_{uv})>a_k\right)<
e^{-\frac{a_k^2}{2q_km_k}+\frac{a_k^3}{2q_k^2m_k^2}},$$
where $$p_k\le q_k=\sum_{(u,v)\in H_k}\frac{p_{uv}}{m_k} <p_{k+1}.$$
We estimate the exponent in case $m_k=n^2$: 
$${-\frac{a_k^2}{2q_km_k}+\frac{a_k^3}{2q_k^2m_k^2}}\le-\frac{\beta^6}{200}\left(\frac{\sqrt{e}^2}{e}\right)^k\frac{Kn^3}{em_k}+
\frac{\beta^9}{2000}\left(\frac{\sqrt{e}^3}{e^2}\right)^k\frac{eKn^5}{m_k^2}<-3n,$$
where we used the definition of $K.$ For $m_k$ being much less than $n^2$ direct substitution gives a
useless bound. For this case we have the useful inequality 

$$\frac{1}{2}{\rm Pr}\left(\sum_{i=1}^{m_k}Y_i>a_k\right) \le 
{\rm Pr}\left(\sum_{i=1}^{n^2}Y_i > \frac{a_k}{2}\right),$$ where 
${\rm Pr}(Y_i=1-q_k)=q_k$ and ${\rm Pr}(Y_i=-q_k)=1-q_k.$ 
This implies that the exponent is at most $-3n$ even in case $m_k <n^2.$

Indeed, let $A, B$ and $C$ be the events that $\sum_{i=1}^{m_k}Y_i >a_k$,
$\sum_{i=1}^{n^2}Y_i > a_k/2$ and $\sum_{i=m_{k+1}}^{n^2}Y_i < -a_k/2$,
respectively. Clearly $A$ and $C$ are independent, and $A \cap \overline{C} \subset B$.
So we have ${\rm Pr}(B) \geq {\rm Pr}(A \cap \overline{C})={\rm Pr}(A){\rm Pr}(\overline{C})$, 
that is ${\rm Pr}(A) \leq {\rm Pr}(B)/{\rm Pr}(\overline{C}) < {\rm Pr}(B)/2$,  
since by $A.1.13$  
$${\rm Pr}\left(\sum_{i=m_{k+1}}^{n^2} Y_i<-\frac{a_k}{2}\right)< e^{-\frac{a_k^2}{8q_k(n^2-m_k)}}<
\frac{1}{2}.$$
With this we have proved that the sum of the weights of the edges of $H_k$ will not be much larger than their expectation with high probability.

Now we estimate the probability that the sum of the weights is much less than their expectation.  
Let us use $A.1.13$ again directly to the sums over $H_k$'s:
$${\rm Pr}\left(\sum_{(u,v)\in H_k}(X_{uv}-{\mathbb E}X_{uv})<-a_k\right)<e^{-\frac{a_k^2}{2q_km_k}}.$$
The exponent in the inequality can be estimated very similarly as before:
$$-\frac{a_k^2}{2q_km_k}\le -\frac{\beta^6}{200}\left(\frac{\sqrt{e}^2}{e}\right)^k\frac{Kn^3}{em_k} <-3n,$$ 
moreover, this bound applies for an arbitrary
$m_k.$

Putting these together we have that
$${\rm Pr}\left(\left|\sum_{(u,v)\in H_k}(X_{uv}-{\mathbb E}X_{uv})\right|>a_k\right)< 2^{-3n}.$$
This implies that

$$\left|\sum_{(u,v)\in H_k}\frac{X_{uv}-{\mathbb E}X_{uv}}{p_{uv}|A||B|}\right|\le \left|\sum_{(u,v)\in H_k}\frac{X_{uv}-{\mathbb E}X_{uv}}{p_{k-1}|A||B|}\right|
\le \frac{a_k}{p_{k-1}|A||B|}\le \frac{\beta}{10}\left(\frac{1}{\sqrt{e}}\right)^k,$$
where the last two inequalities hold with probability at least $1-2^{-3n}$ for a given pair of sets $A$ and $B$ if $n$ is sufficiently large.
Since 
$$\left| \sum_{u \in A, v\in B}\frac{X_{uv}}{p_{uv}|A||B|} \right| =\left| \sum_{k=1}^{\log n} \sum_{(u,v)\in H_k}\frac{X_{uv}}{p_{uv}|A||B|}  \right| \le  \sum_{k=1}^{\log n} 
\left|\sum_{(u,v)\in H_k}\frac{X_{uv}}{p_{k-1}|A||B|}\right|$$
and $$\sum_{k=1}^{\log n}\frac{1}{10}\left(\frac{1}{\sqrt{e}}\right)^k\le \frac{1}{2},$$
the claimed bound follows with high probability. 
\hfill \qedf

\bigskip

\noindent {\bf Remark.} It is very similar to prove that with high probability $|\sum_{i, j} \rho_G(i,j) -\binom{n}{2}|=o(n),$ 
we omit the details. From this it follows that rescaling the above edge weights by a factor of $(1+o(1))$ and letting $\mu \equiv 1$ 
provides us $\beta$-quasi-random weights for most graphs from $G(n, p_{ij})$ such that $\mu(V)=n$ and 
$\rho_G(V,V)=2 \binom{n}{2}.$ That is, with high probability we can apply the Regularity Lemma for any $F \subset G$, where 
$G \in G(n, p_{ij}).$

\subsection{Simple examples for defining vertex and edge weights}

When defining the notion of weighted quasi-randomness and weighted regularity, we mentioned that choosing 
$\mu \equiv 1$ and $\rho \equiv 1/\delta_G$ gives back the old definitions of quasi-randomness and regularity. 
In the previous section we saw an example when we needed different edge weights, but $\mu$ was identically one. 

Let us consider a simple example in which $\mu$ has to take more than one value. Let $G$ be a star on $n$ vertices, that is, the vertex $v_1$ is adjacent to the vertices 
$v_2, \ldots, v_n,$ and $v_i$ has degree 1 for $i \ge 2.$ We let $\mu(v_1)=1/2$ and $\mu(v_i)=1/(2(n-1))$ for $i \ge 2,$ and choose $\rho_G \equiv n/2.$ With these choices $G$ 
is easily seen to be quasi-random, moreover, it is weighted regular.

A more sophisticated example relates weighted regularity with volume regularity, the latter introduced by Alon et al~\cite{A}. Given a graph $G,$ 
the volume of a set $S \subset V(G)$ is 
defined as ${\rm vol}(S)=\sum_{v \in S}{\rm deg}(v).$ For disjoint sets $A, B \subset V(G),$
the $(A,B)$ pair is $\varepsilon$-volume regular if $\forall X \subset A, Y \subset B$ satisfying ${\rm vol}(X) \ge \varepsilon \cdot {\rm vol}(A)$ and ${\rm vol}(Y) \ge \varepsilon \cdot {\rm vol}(B)$
we have
$$\left |e(X,Y)-\frac{e(A,B)}{{\rm vol}(A){\rm vol}(B)}{\rm vol}(X){\rm vol}(Y)\right| < \varepsilon \cdot {\rm vol}(A){\rm vol}(B)/{\rm vol}(V).$$
In order to show that volume regularity is a special case of weighted regularity, we may choose edge and vertex weights as follows:
for all $A, B \subset V$ 
let $$\rho (A, B)=e(A, B)\frac{n^2}{{\rm vol}(V)}$$ and 
$$\mu(A)=\frac{n \cdot {\rm vol}(A)}{{\rm vol}(V)}.$$
It is easy to see that if the $(A,B)$ pair is $(\mu, \rho)$-weighted $\varepsilon$-regular with the above choice for $\mu$ and $\rho$, 
then it is also $\varepsilon$-volume regular.

We will consider a method for constructing vertex weights in the next section, that will also demonstrate the difference between volume regularity 
and weighted regularity.

\subsection{Natural weighting of $K_n$}

Assume that $|V|=n.$ Let the vertex weight function $\mu: V \rightarrow \mathbb R^+$ be defined such that 
$\mu(V)=n.$
Then we define the {\it natural weighting of the edges of $K_n$ with respect to $\mu$} as follows:  we let 
$\rho(u,v)=\mu(u)\cdot \mu(v)$ for all $u, v \in V,$ $u\neq v.$
We show that these weight functions determine a quasi-random weighting of $K_n.$ Let $A, B \subset V$ such that $A \cap B =\emptyset.$ Then
$$\frac{\rho(A, B)}{\mu(A)\mu(B)}=\frac{\sum_{u\in A}\sum_{v\in B}\mu(u)\mu(v)}{\mu(A)\mu(B)}=\frac{\mu(A)\mu(B)}{\mu(A)\mu(B)}=1,$$ independent
of the weights of $A$ and $B.$ Since we allow loop edges, $\rho(V, V)=\mu(V)\mu(V).$ 
Recalling the definition of quasi-randomness it is easy to see that the natural weighting of $K_n$ is {\em always} quasi-random. 

Note, that natural weighting resembles to Definition~\ref{quasirandom}, where the lower bounds on $\mu(A)$ and
$\mu(B)$ are dropped. It is closely related to the random model $G(\mathbf{w})$, see e.g. in~\cite{CL}. Here  
$\mathbf{w}=(w_1, \dots, w_n)$ is the expected degree sequence of $G(\mathbf{w})$ with vertex set $\{1, 2, \ldots, n\}.$ The edges of $G(\mathbf{w})$ 
are drawn independently, and the probability of including the edge $ij$ is $w_i w_j/\sum_iw_i$. Of course, the model 
$G(\mathbf{w})$ is the special case of $G(n, p_{i j})$, and Lemma~\ref{pij} holds {\em without} any conditions.

Let $u$ be an arbitrary vertex and $A\subset V.$ Then the {\it weighted degree of $u$ into $A$} 
in the graph $F \subset K_n$ is defined to be 
$$dw_F(u, A)=\sum_{v\in A}{\mathbf 1}_F(u, v)\mu(v)=\mu(N_F(u,A)),$$ where $N_F(u, A)$ denotes the neighborhood of $u$ 
in the set $A.$ In particular the weighted degree of $u$ in $F$ is $$dw_F(u)=\sum_{v \in V}{\mathbf 1}_F(u, v)\mu(v)=\mu(N_F(u)).$$ 
We also have that 
$$\rho_F(A, B)=\sum_{u \in A}\sum_{v \in B}{\mathbf 1}_F(u,v)\rho(u, v)=\sum_{u \in A}dw_F(u, B)$$ and 
$$\rho_F(V, V)=\sum_{u \in V}dw_F(u).$$
We define the weighted density of a weighted $\varepsilon$-regular $(A,B)$ pair to be $$\frac{\rho_F(A, B)}{\mu(A)\mu(B)}.$$ 
We have the following lemma.

\begin{lemma}\label{fok}
Let $(A,B)$ be a weighted $\varepsilon$-regular pair relative to the natural weighting of $K_n$ with weighted density 
$\gamma \gg \varepsilon.$ Let $A'\subset A$ contain only such vertices that have weighted degree less than 
$(\gamma-\varepsilon)\mu(B)$ in the pair. Then $\mu(A')<\varepsilon \mu(A).$ 
\end{lemma}

\noindent {\bf Proof:} Assume on the contrary that the set of `low-degree' vertices has a large weight. 
Observe that $\varepsilon$-regularity implies that $$\frac{\rho_F(A',B)}{\mu(A')\mu(B)} >\gamma-\varepsilon$$ if $\mu(A')>\varepsilon \mu(A).$
Using our assumption we get the following:
$$\gamma-\varepsilon<\frac{\rho_F(A',B)}{\mu(A')\mu(B)} =\frac{\sum_{u \in A'}\mu(u)dw_F(u, B)}{\mu(A')\mu(B)}<\frac{\sum_{u \in A'}\mu(u)(\gamma-\varepsilon)\mu(B)}{\mu(A')\mu(B)}=
\gamma-\varepsilon,$$ which is clearly a contradiction.
\hfill\qedf

\medskip

Let $A, B_1, B_2, \ldots, B_k$ be disjoint subsets of $V(F),$ and assume that $(A, B_i)$ is a weighted $\varepsilon$-regular pair relative to a natural weighting of $K_n$ 
with weighted density at least $\gamma$ for every $i.$ 
Set $\delta=\gamma-\varepsilon.$ Let $0<s$ be an integer constant, and assume that $\delta^s \gg \varepsilon.$ 

\begin{lemma}\label{nagyfok}  Assume that $A'\subset A$ with $\mu(A') > 2k\varepsilon \mu(A).$ Then there exist vertices $u_1, u_2, \ldots, u_s \in A' $ such that 
$$\mu(\cap_{1\le i \le s} N_F(u_i, B_j)) \ge \delta^s\mu(B_j)$$ for every $1 \le j \le k.$ 
\end{lemma}

\noindent{\bf Proof:} We find the $u_i$ vertices one by one. For $u_1$ we have that the weight of vertices of $A$ with weighted degree at most $\delta \mu(B_1)$ is at most
$\varepsilon \mu(A)$ using Lemma~\ref{fok}. Discard these low-degree vertices from $A'$, then use the regularity condition again, this time for $B_2.$ 
We find that the weight of vertices having small degree
into $B_1$ or $B_2$ is at most $2\varepsilon \mu(A).$ Iterating this procedure we get that the weight of vertices that do not have large degree into 
at least one $B_i$ set is at most $k\varepsilon\mu(A)< \mu(A').$
Pick any of the large degree vertices from $A',$ this is our choice for $u_1.$ 

Next we repeat the process for finding $u_2,$ with the difference that we look for a vertex that have large degree into the sets $N_F(u_1, B_j)$ for every $j.$ Since 
$\mu(N_F(u_1, B_j))\ge \delta \mu(B_j) \gg \varepsilon \mu(B_j),$ the same procedure will work. Applying Lemma~\ref{fok} we can find many vertices in $A'-u_1$ such that the weighted degree 
of all of them into $B_j$ is 
at least $\delta \mu(N_F(u_1, B_j))\ge \delta^2\mu(B_j)$ for every $j.$ Pick any of these, this vertex is $u_2.$   

When it comes to finding $u_i$ we will work with the sets $A'-\{u_1, \ldots, u_{i-1}\}$
and $\cap _{t\le i-1}N_F(u_t, B_j)$ for $1 \le j \le k.$ Using induction it is easy to show that $$\mu(\cap _{t\le i-1}N_F(u_t, B_j))\ge \delta^{i-1}\mu(B_j)$$ for every $j.$ 
Since $\delta^{s} \gg \varepsilon,$ we can iterate this procedure until we find all the vertices $u_1, \ldots, u_s$. 
\hfill \qedf

\medskip

Assume now that there are $q$ clusters, $W_1, W_2, \ldots, W_q \subset V(F)$ such that $\mu(W_i)=m+o(m)$ for all $i$ (here $m=constant \cdot n$) and all  the $(W_i, W_j)$ pairs are 
weighted $\varepsilon$-regular relative to a natural weighting of $K_n$
with density at least $\gamma.$ That is, we have a super-clique $Cl_q$ on $q$ clusters. 

Next we construct the graph $K^s_q,$ a blown-up clique as follows. First, we have $q$ disjoint $s$-element set of vertices, this is the vertex set of $K^s_q.$ Then we connect any two vertices
if they belong to different vertex sets. Before we state an embedding result, we need a simple lemma, the proof is left for the reader.

\begin{lemma}\label{slicing}
Let $(A,B)$ be a weighted $\varepsilon$-regular pair with density $d,$ and for some $\alpha$ let $A' \subset A$ with $\mu(A')\ge \alpha\mu(A)$ and $B' \subset B$ 
with $\mu(B')\ge \alpha\mu(B).$ Then $(A',B')$ is a weighted $\varepsilon'$-regular pair with $\varepsilon'=\max\{\varepsilon/\alpha, 2\varepsilon\}$ and density $d'\ge d-\varepsilon.$
\end{lemma}

We have the following embedding lemma:

\begin{lemma}\label{kqs}
Let $\delta=\gamma - 2\varepsilon.$ If $\delta^{qs}\gg \varepsilon$ then $K^s_q\subset Cl_q.$
\end{lemma}

\noindent {\bf Proof:} First, apply Lemma~\ref{nagyfok} with $A=W_1$ and $B_j=W_{j+1}$ for $1 \le j \le q-1.$ We find the the vertices $u_1^1, u_2^1, \ldots, u_s^1 \in W_1$ such that
$$\mu(\cap_{1 \le i \le s}N_F(u_i^1, W_j)) \ge \delta^s\mu(W_j).$$ Let $W_j^2=\cap_{i\ge 1}N_F(u_i^1, W_j)$, then $\mu(W_j^2)\ge \delta^s\mu(W_j)\gg \varepsilon \mu(W_j)$
for every $j \ge 2.$ 

Next let $A=W_2^2$ and $B_j=W^2_{j+2}$ for $1 \le j \le q-2.$ Using Lemma~\ref{slicing} we have that the new $(A, B_j)$ pairs are all
weighted $\varepsilon/\delta^s$-regular with density at least $\gamma-\varepsilon.$ Hence, we can apply Lemma~\ref{nagyfok} again, and find $u_1^2, u_2^2, \ldots, u_s^2 \in W_2^2$
such that $$\mu(\cap_{1 \le i \le s}N_F(u_i^2, W_j^2)) \ge \delta^s\mu(W_j^2)\ge \delta^{2s}\mu(W_j) \gg \varepsilon \mu(W_j)$$ for $3 \le j \le q.$

Continuing this process, in the $k$th step we will work with the $W_j^{k}$ sets when applying Lemma~\ref{nagyfok}. These sets are defined recursively as follows: 
$W_j^k=\cap_{i\ge 1}N_F(u_i^{k-1}, W_j^{k-1})$
and $\mu(W_j^k) \ge \delta^{(k-1)s}\mu(W_j)$ for every $k+1 \le j \le q.$ Moreover, the $(W_k^k, W_j^k)$ pairs will be $\varepsilon/\delta^{(k-1)s}$-regular with density at least 
$\gamma-\varepsilon$ for every $k+1 \le j \le q.$
 
In the last step, when $k=q-1,$ there are only two sub-clusters left, $W_{q-1}^{q-1}$ and $W_{q}^{q-1}.$ The pair $(W_{q-1}^{q-1},W_{q}^{q-1})$ will be weighted 
$\varepsilon/\delta^{(q-2)s}$-regular
with density at least $\gamma -\varepsilon.$ It is easy to find a $K_{s,s}$ (a complete bipartite graph) in this regular pair using Lemma~\ref{nagyfok}. Clearly, we constructed
the desired $K_q^s$ graph. \hfill \qedf

\subsection{A weighted version of the Erd\H{o}s-Stone theorem}

Let $t_{q-1}(n)$ be the number of edges in the Tur\'an graph $T_{n, q-1}$ on $n$ vertices. That is, $T_{n, q-1}$  has the largest number of edges such that it does not 
contain a $K_q.$ It is well known that
$$\lim_{n \rightarrow \infty} \frac{t_{q-1}(n)}{{n \choose 2}}=\frac{q-2}{q-1}.$$ 

The Erd\H{o}s-Stone theorem states that if one has at least $t_{q-1}(n)+\gamma n^2$ edges (where $\gamma >0$ is a constant) in a graph $F$ on $n$ vertices then $F$ has a 
$K_q^s$ for any given natural number $s.$ In this section we show a weighted version. We take a natural weighting of $K_n,$ and prove that 
if the total edge weight in $F \subset K_n$ is large then $F$ has a large blown-up clique. We remark that there are other results in the literature on the extremal theory of weighted
graphs, see e.g. ~\cite{BT} by Bondy and Tuza and~\cite{FK} by F\"uredi and K\"undgen, although the setup of these papers is different from ours. 

\begin{theorem}\label{ES}
For all integers $q\ge 2$ and $s\ge 1$ and every $\gamma>0$ there exists an integer $n_0$ such that the following holds. Take the natural weighting of $K_n$ with vertex weight function $\mu$ 
and assume that $\mu(V)=n\ge n_0.$ Let $F \subset K_n. $ If the total edge weight of $F$ is at least $t_{q-1}(n)+\gamma n^2$ then $F$ contains $K_q^s$ as a subgraph.  
\end{theorem}

\noindent {\bf Proof:}
We begin with applying the weighted Regularity Lemma with parameters $\varepsilon \ll \min\{(\gamma-\varepsilon)^{qs}, 1/s, 1/q\}$ and $L\gg 1/\varepsilon.$ 
We get an $\varepsilon$-regular partition with 
clusters $W_0, W_1, \ldots, W_{\ell}.$ Let us construct the {\it reduced graph $F_r$} as follows. The vertices of $F_r$ are identified by the $\ell$ non-exceptional clusters. 
We have an edge between two vertices of $F_r$ if the corresponding two clusters give an $\varepsilon$-regular pair with density at least $\gamma.$ Hence, when we construct $F_r$
we lose edges of $F$ as follows: (1) edges that are incident with some vertex of $W_0,$ (2) edges that connect two vertices that belong to the same non-exceptional cluster, 
(3) edges that
are in some irregular pair, (4) edges that are in regular pairs with small density. 

The outline of the proof is as follows. We will show that the loss in edge weight is small, hence, $F_r$ will have many edges.  By Tur\'an's Theorem we will have a $q$-clique
in $F_r.$ Then we apply Lemma~\ref{kqs}, and conclude the existence of a $K_q^s$ in $F.$

(1) The total weight of edges that are incident with some vertex of $W_0$ can be estimated as follows: 
$$\rho_F(W_0, V) \le \rho(W_0,V) =\sum_{w\in W_0}\sum_{v\in V}\rho(w,v)\le \mu(W_0)\mu(V)\le \varepsilon n^2.$$

(2) The non-exceptional clusters have weight $(n-\varepsilon n)(1+o(1))/\ell.$ The total weight of edges inside non-exceptional clusters is at most 
$$\frac{1}{2}\sum_{1 \le i \le \ell}\sum_{u \in W_i}\sum_{v \in W_i-u}\rho(u,v)=\frac{1}{2}\sum_{1 \le i \le \ell}\sum_{u \in W_i}\sum_{v \in W_i-u}\mu(u)\mu(v)\le$$
$$\frac{1}{2}\sum_{1 \le i \le \ell}\mu(W_i)^2=\frac{n^2}{\ell}(1+o(1)).$$
Since $\ell \ge L \gg 1/\varepsilon,$ we have that the total edge weight inside non-exceptional clusters is less than $\varepsilon n^2.$

(3) Assume that $(W_i, W_j)$ is an irregular pair. Then $$\rho_F(W_i, W_j)\le \sum_{u \in W_i}\sum_{v \in W_j}\rho(u, v)=
\sum_{u \in W_i}\sum_{v \in W_j}\mu(u) \mu(v)=\mu(W_i)\mu(W_j)=\frac{n^2}{\ell^2}(1+o(1)).$$ Since the number of irregular pairs is at most $\varepsilon \ell^2,$
we get that the total edge weight in irregular pairs is at most $$\varepsilon \ell^2 \frac{n^2}{\ell^2}(1+o(1))<2\varepsilon n^2.$$

(4)  If the density of an $\varepsilon$-regular pair $(W_i, W_j)$ is small then we have the following inequality:
$$\rho_F(W_i, W_j)\le \mu(W_i)\mu(W_j)\gamma=\frac{n^2}{\ell^2}(1+o(1))\gamma.$$ Since there can be at most ${\ell \choose 2}$ pairs, the total edge weight in 
low density pairs is less than $2\gamma n^2/3.$

Putting together, we get that the total weight of edges that we lose when applying the weighted Regularity Lemma is at most $(4\varepsilon +2\gamma/3) n^2<
3\gamma n^2/4.$ Hence, the total edge weight in the high-density regular pairs of $F_r$ is at least $t_{q-1}(n)+\gamma n^2/4.$
The total weight of edges in a regular pair is $(1+o(1))n^2/\ell^2.$ Assume that $e(F_r) \le \frac{q-2}{q-1}\frac{\ell^2}{2},$ then the total edge weight would be at 
most $\frac{q-2}{q-1}\frac{n^2}{2}(1+o(1)).$ Since we have a larger edge weight in what is left after applying the Regularity Lemma, using Tur\'an's theorem, we get that
$F_r$ contains a $K_q.$ Every pair in this clique is a high-density $\varepsilon$-regular pair, hence, we can apply Lemma~\ref{kqs} and find the blown-up clique. 
\hfill \qedf

\medskip

\noindent{\bf Remarks.}
One can arrive at the same conclusion perturbing the edge weights a little. Let $D>1$ be a fixed constant. Multiply the weight of the edge $e$ by any number 
$c_e \in [1/D, D].$ The resulting weighted graph will be quasi-random, and it is an easy exercise to show that one still have Theorem~\ref{ES}.   

One can also show the weighted version of the Erd\H{o}s-Stone-Simonovits theorem, a stability version of the above. Let $\mathcal H$ be a family of forbidden subgraphs
having chromatic number $q.$ 
Assume that the total edge weight in $F$ is close to $t_{q-1}(n),$ but $F$ does not contain some graph $H \in {\mathcal H}.$ Then $F_r,$ the reduced graph
cannot have a clique on $q$ vertices, but the number of edges in it will be close to $t_{q-1}(\ell).$ This implies that $F_r$ is close to a Tur\'an graph $T_{\ell, q-1},$
and that in turn implies that the vertex set of $F$ can be partitioned into $q-1$ disjoint vertex classes in the following way: the vertex classes all have weight $\approx n/(q-1),$
the total weight of edges inside vertex classes is very small, and the weighted density of edges for every pair of disjoint classes  is close to one.

\subsection{Emphasized sets}\label{fontos}
One cannot avoid to have an exceptional cluster $W_0$ when applying the Regularity Lemma. That is, a linear number of vertices could be 
discarded in certain cases, a well-know example is the so called half-graph. In general we don't have a control on what is put into the exceptional 
cluster. However, using vertex weights one can at least partly determine the set of discarded vertices. In what follows we show how to the natural weighting
of $K_n$ in order to have that the majority of some given emphasized set is put into non-exceptional clusters after applying the weighted Regularity Lemma, 
even if the set is of size $o(n).$ In fact we will do it for several emphasized sets at the same time. Notice that using the volume regularity concept of~\cite{A} one
may discard all vertices with small degrees.

Assume that $k$ is a fixed constant, $V$ is partitioned into the disjoint sets $S_1, S_2, \ldots, S_k,$ and let $n=|V|.$ Further assume 
that $s_i \rightarrow \infty$ as $n  \rightarrow \infty.$ Let $s_i=|S_i|$ for every $i.$ Define the following weighting of the vertices of $V$: for $v\in S_i$ we let $$\mu(v)=\mu_i=\frac{n} {ks_i}.$$
Observe that $$\sum_{v \in S_i}\mu(v)=\frac{n}{k},$$ thus, the total weight of the vertices is $n.$
Let $v \in S_i$ and $w\in S_j.$ The weight of the pair $(v,w)$ is $$\rho(v,w)=\rho_{ij}=\mu(v)\mu(w)=\frac{n^2}{k^2s_is_j}.$$ We showed above that $K_n$ equipped with such
vertex and edge weights is a quasi-random graph. We call this weighting the {\it natural weighting of $K_n$ with emphasized sets $S_1, S_2, \ldots, S_k.$}
 
We can apply Theorem~\ref{grl} for some $F$ relative to the natural weighting of $K_n.$  Choose $\varepsilon$ so that $k \ll 1/\varepsilon.$ 
Since $\mu(W_0)\le \varepsilon n \ll n/k,$ we get that for all $i$ the majority of the vertices of $S_i$ are in non-exceptional clusters. 

We remark that it is possible to define vertex weights not only for $G=K_n,$ but for much sparser quasi-random graphs when emphasizing subsets of $V.$ 
For example, assume that $G\in G(n,p_{ij}),$ and $V$ is partitioned into the disjoint sets $S_1, S_2, \ldots, S_k.$ Then one will have the vertex weights of the above example, 
but the edge weights will be different: $$\rho(v_i,v_j)=\mu(v_i)\mu(v_j)\frac{1}{p_{ij}}=\frac{n^2}{k^2s_qs_tp_{ij}}$$ whenever $v_i \in S_q$ and $v_j \in S_t.$ We leave the details for the reader.

\medskip

\noindent {\bf Acknowledgment}
The authors thank P\'eter Hajnal and Endre Szemer\'edi  for the helpful discussions.
We also thank the anonymous referees for the numerous suggestions which improved 
the presentation of the paper.

\end{document}